\documentclass[11pt]{article}
\usepackage[margin=0.75in]{geometry}
\usepackage{bookmark}
\usepackage{hyperref}
\usepackage{graphicx}
\usepackage{float}
\usepackage{amsfonts}
\usepackage{amsmath,amssymb}
\usepackage{epstopdf}
\usepackage{subfigure}
\usepackage{algorithm, algorithmic}
\usepackage[font=footnotesize]{caption}
\usepackage{color}
\usepackage{wrapfig}
\usepackage{marvosym}
\usepackage{caption}
\usepackage{cleveref}
\usepackage{ragged2e}
\usepackage[scaled]{helvet} 
\usepackage{fourier}

\newcommand{\C}{\mathbb{C}}

\newcommand{\cS}{\mathcal{S}}
\newcommand{\cP}{\mathcal{P}}
\newcommand{\cR}{\mathcal{R}}

\newcommand{\bD}{\boldsymbol{D}}

\newcommand{\br}{\boldsymbol{r}}
\newcommand{\bx}{\boldsymbol{x}}

\newcommand{\balpha}{\boldsymbol{\alpha}}
\newcommand{\bA}{\boldsymbol{A}}
\newcommand{\by}{\boldsymbol{y}}

\newcommand{\be}{\boldsymbol{e}}

\newcommand{\bz}{\boldsymbol{z}}
\newcommand{\bw}{\boldsymbol{w}}
\newcommand{\bv}{\boldsymbol{v}}

\newcommand{\bal}{\boldsymbol{\alpha}}

\newcommand{\bPhi}{\boldsymbol{\Phi}}
\newcommand{\bI}{\boldsymbol{I}}

\newcommand{\supp}{\mathrm{supp}}
\newcommand{\argmin}{\mathop{\mathrm{argmin}}}

\numberwithin{equation}{subsection}

\title{Practical approximate projection schemes\\ in greedy signal space methods}

\author{Chris Garnatz\thanks{Pomona College, Claremont CA 91711}, Xiaoyi Gu\thanks{Univ. of California, Los Angeles CA 90095}, Alison Kingman\thanks{Harvey Mudd College, Claremont CA 91711}, James LaManna\thanks{Pitzer College, Claremont CA 91711}, Deanna Needell\thanks{Claremont McKenna College, Claremont CA 91711{\tt\small dneedell@cmc.edu}}, Shenyinying Tu\thanks{Univ. of California, Los Angeles CA 90095} 
\thanks{This work was supported by NSF grant DMS-$1045536$, NSF CAREER grant $\#1348721$, and the Alfred P. Sloan Fellowship.}%
}

\date{\today}

\begin{document}
\maketitle

\begin{abstract}
Compressive sensing (CS) is a new signal acquisition paradigm which shows that far fewer samples are required to reconstruct sparse signals than previously thought.  Although most of the literature focuses on signals sparse in a fixed orthonormal basis, recently the Signal Space CoSaMP (SSCoSaMP) greedy method was developed for the reconstruction of signals compressible in arbitrary redundant dictionaries.  The algorithm itself needs access to approximate sparse projection schemes, which have been difficult to obtain and analyze.  This paper investigates the use of several different projection schemes and catalogs for what types of signals each scheme can successfully be utilized.  In addition, we present novel hybrid projection methods which outperform all other schemes on a wide variety of signal classes.
\end{abstract}

\section{Introduction}
Compressive Sensing (CS) is an emerging field that seeks to optimize signal acquisition and reconstruction. 
CS relies on the fact that most real world signals are themselves sparse or compressible with respect to some fixed basis.  Traditionally, this basis is assumed to be orthonormal.  However, most signals in practice are instead sparse with respect to non-orthonormal, highly overcomplete dictionaries.  CS has only recently begun to develop technologies for this framework, and we investigate some of them in this paper.

We will use the following notation to frame the fundamental problem CS seeks to answer. The signal in question is $\bx \in \C^n$. The measurement vector $\by \in \C^m$ is related to the signal by the measurement matrix $\bA \in \C^{m \times n}$ and additive noise $\be \in \C^m$, where $\by = \bA\bx + \be$ and $m \ll n$. We impose a sparsity condition on $\bx$ by the equation $\bx = \bD \bal$, where $\bal \in \C^d$ with $|\supp(\bal)| \leq k \ll n$ is a $k$-sparse representation of the signal $\bx$ with respect to the dictionary $\bD \in \C^{n \times d}$.  Note that when $d = n$ and $\bD = \bI_n$ is the identity, then we have the case where $\bx$ is itself sparse.

One can recover a sparse vector from its measurements by solving an $\ell_0$-minimization problem which simply searches for the sparsest vector which yields the same measurements.  Although this is an NP-hard problem, the seminal work of Cand\`{e}s, Romberg and Tao \cite{RefWorks:48, RefWorks:47} shows that a sparse vector $\bx$ can be recovered using the relaxed $\ell_1$-minimization method from the measurement vector $\by$ as long as the measurement matrix $\bA$ satisfies a condition known as the Restricted Isometry Property (RIP).  To say $\bA$ satisfies the (RIP) of order $k$ is to say the following: for a constant $\delta_k \in (0,1)$, the condition

\begin{equation}
\label{RIP}
{(1-\delta_k)}||\bx||_2^2 \leq ||\bA \bx||_2^2  \leq {(1+\delta_k)}||\bx||_2^2
\end{equation}

\noindent holds for all $\bx$ with $|\supp(\bx)| \leq k$. Finding a measurement matrix $\bA$ that satisfies the RIP for a sparsity level $k$ is a nontrivial endeavor. Fortunately, if the entries of $\bA$ are drawn from a subgaussian distribution and we use $m\geq Ck/\log (n)$ measurements, then $\bA$ satisfies the RIP with high probability \cite{RefWorks:444,RefWorks:285}.  Similar results hold for subsampled Discrete Fourier Transform (DFT) matrices and others with a fast-multiply \cite{RefWorks:285, rauhut2012restricted}.  Thus a measurement matrix $\bA$ with sufficient properties for signal recovery can be practically constructed.  Greedy approaches such as OMP~\cite{Paper9, RefWorks:150}, ROMP~\cite{RefWorks:44}, CoSaMP~\cite{NeedeT_CoSaMP}, and IHT~\cite{PaperIHT} which identify elements of the signal support iteratively also provide robust reconstruction guarantees similar to those of the $\ell_1$ approach.  We point the reader to these references for details about each of the CS algorithms.

\section{Signal Space Methods}

Existing CS literature has largely addressed the scenarios where $\bx$ is sparse or compressible with respect to an orthonormal dictionary $\bD$.  When $\bD$ is orthonormal, these methods simply recover $\balpha$, which naturally yields $\bx$.  
Most signals in practical applications, however, are compressible instead with respect to a highly overcomplete dictionary such as an oversampled DFT, Gabor frame, or many of the redundant frames used in image processing.
Unfortunately, when $\bD$ is not orthonormal, these CS methods fail both theoretically and empirically. 

In \cite{RefWorks:23,Paper5}, Davenport et. al. propose a new method called Signal Space CoSaMP (SSCoSaMP), specifically designed to recover the signal $\bx$ when $\bD$ is not guaranteed to be orthonormal. They utilize a generalization of the RIP, the $D$-RIP criterion introduced in \cite{RefWorks:60} for the $\ell_1$-analysis approach, to ensure signal recovery.  In contrast to other CS algorithms, SSCoSaMP doesn't return $\hat{\balpha}$ to obtain $\hat{\bx}$, but instead opts for a ``signal-focused" approach and attempts to directly calculate and return $\hat{\bx}$.

The two key steps within SSCoSaMP require finding the best $k$-sparse approximation of a vector $\bz$ with respect to the arbitrary dictionary $\bD$:
$$
\Omega_{\text{opt}} := \argmin_{\Lambda: |\Lambda| = k} \|{\bz - \cP_{\Lambda} \bz}\|_2,
$$
where $\cP_{\Lambda}$ denotes the projection onto the span of the columns of $\bD$ indexed by $\Lambda$.  Unfortunately, finding the support of such an optimal sparse projection vector is an NP-hard problem itself. Instead, in the interest of computational efficiency, SSCoSaMP computes a \emph{near-optimal} approximation~\cite{Paper5,Paper6} of the support via a simpler CS algorithm, which is noted by the function $\cS_{\bD}$.  Stunningly, SSCoSaMP is able to accurately recover signals with great success even when $\bD$ is very overcomplete and redundant, using a standard CS method like OMP, CoSaMP or $\ell_1$-minimization for the near-optimal projection $\cS_{\bD}$. The SSCoSaMP algorithm is described in Algorithm~\ref{ccosamp}, see~\cite{Paper5} for a detailed description.

\begin{algorithm}[H]
\caption{Signal-Space CoSaMP (SSCoSaMP)}\label{ccosamp}
\begin{algorithmic}

\STATE \textbf{Input:} $\bA$, $\bD$, $\by$, $k$, stopping criterion 
\STATE \textbf{Initialize:} $\br = \by$, $\bx_0 = 0$, $\ell = 0$, $\Gamma = \varnothing$

\WHILE{not converged}
\STATE
\begin{tabular}{ll}

\textbf{Proxy:} & $\widetilde{\bv} = \bA^*\br$  \\ 
\textbf{Identify:} & $\Omega = \cS_{\bD}(\widetilde{\bv},2k)$ \\
\textbf{Merge:} & $T = \Omega \cup \Gamma$ \\
\textbf{Update:} & $\widetilde{\bw} = \argmin_{z}||\by - \bA \bz||_2 \quad   \mathrm{s.t.}  \quad  \bz \in \cR(\bD_T)$	\\
\textbf{Prune:} & $\Gamma = \cS_{\bD}(\widetilde{\bw}, k)$\\
 & $\bx_{\ell+1} = \cP_{\Gamma}\widetilde{\bx}$\\
& $\br = \by - \bA\bx_{\ell + 1}$\\
& $\ell = \ell + 1$\\

\end{tabular}
\ENDWHILE

\STATE \textbf{Output:} $\hat{\bx} = \bx_{\ell}$

\end{algorithmic}
\end{algorithm}

\begin{figure}

\centering
\hspace{-.075\linewidth}
\begin{minipage}{.45\linewidth}
\centering
\begin{tabular}{cc}
\raisebox{25mm}{ {(a)}} &\includegraphics[width=3in]{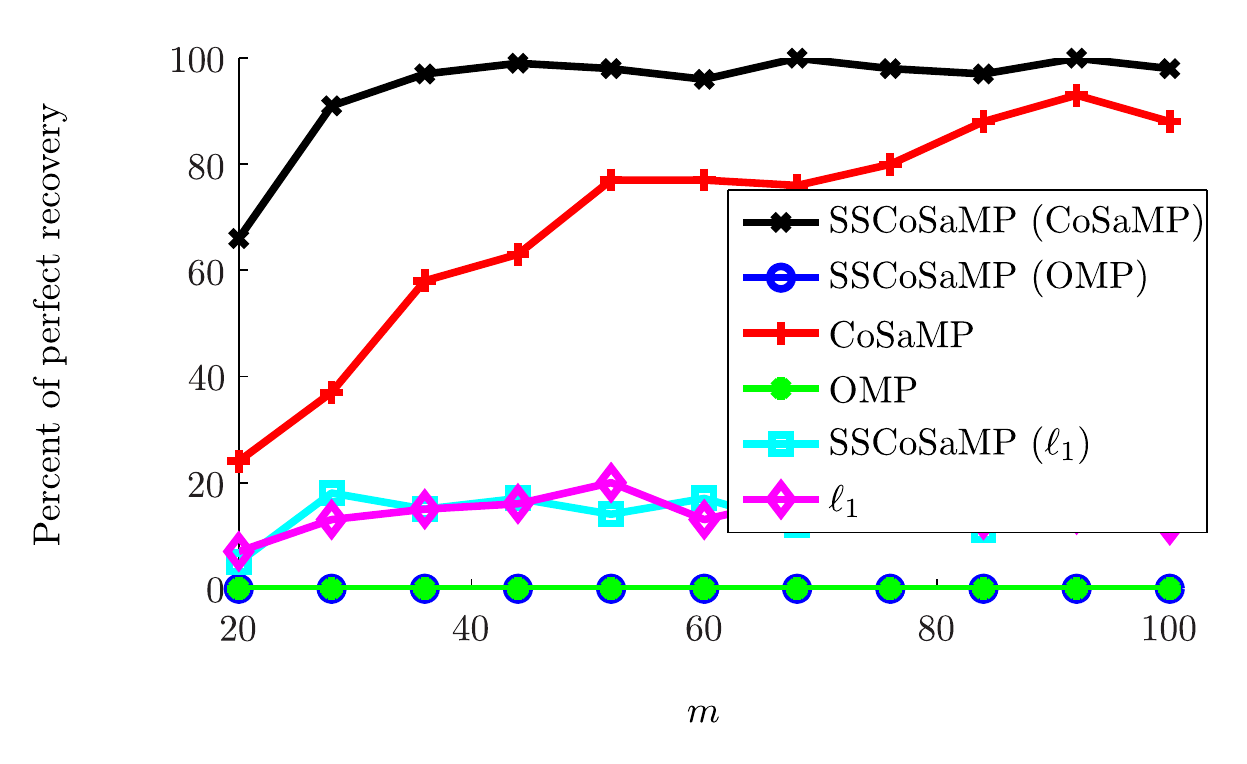}
\end{tabular}
\end{minipage}
\hspace{.05\linewidth}
\begin{minipage}{.45\linewidth}
\centering
\begin{tabular}{cc}
\raisebox{25mm}{ {(b)}} &\includegraphics[width=3in]{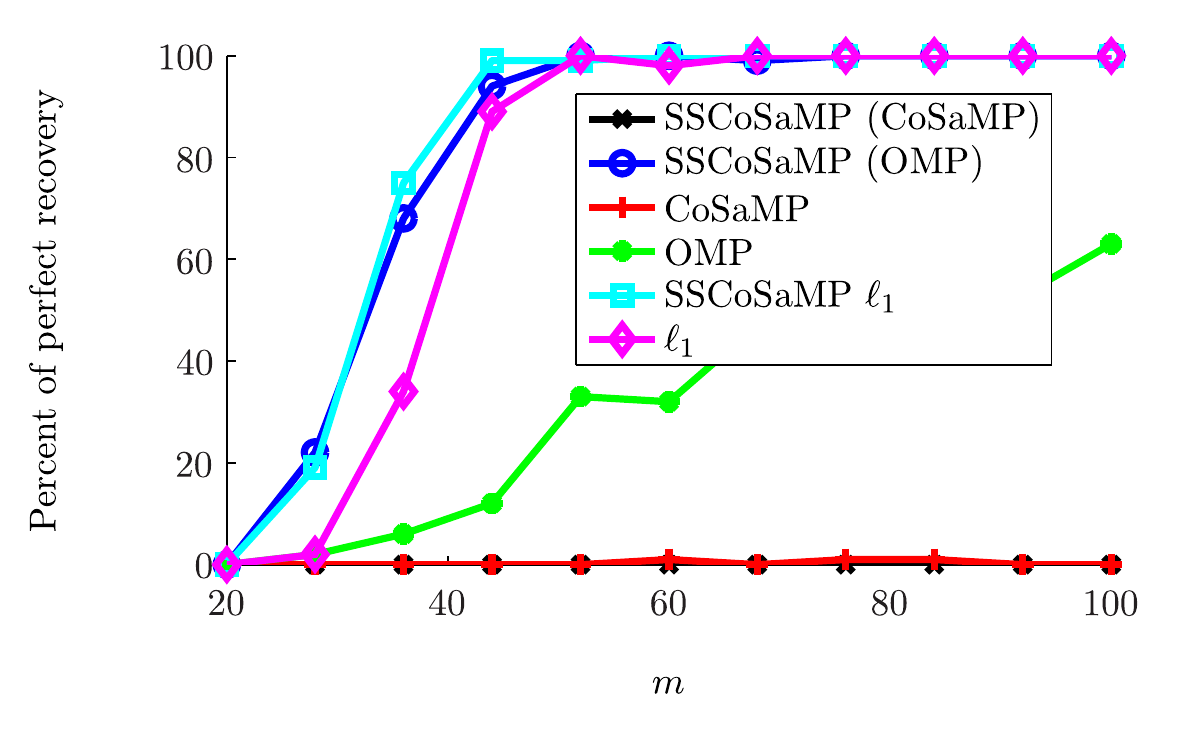}
\end{tabular}
\end{minipage}
   \caption{ Percent of perfect recovery for different SSCoSaMP variants when the nonzero entries in $\balpha$ are clustered together (left) and when the nonzero entries in $\balpha$ are well-separated (right).
   \label{sepvclust}}
\end{figure}

SSCoSaMP is a state of the art algorithm, but in its present form, its recovery performance varies greatly with what CS algorithm is used for the near-optimal projection and with the structure of the sparse coefficient vector $\balpha$. This phenomenon is highly evident in the simulation results presented by the authors of \cite{Paper5}, which we recreate in Figure \ref{sepvclust}, which shows percent of perfect recovery as a function of the number of measurements $m$.\footnote{We define perfect recovery as signal to noise ratio $20\log(\|\bx\|_2/\|\bx-\hat{\bx}\|_2)$ greater than 100db.} Here, the dictionary $\bD \in \mathbb{C}^{n\times d}$ is a $4\times$ overcomplete DFT and $\bA \in\mathbb{R}^{m\times n}$ has standard normal entries.  For all $100$ trials, $k = 8$, $n = 256$, and $d = 1024$.  Here and throughout this paper, the algorithm in parentheses following ``SSCoSaMP" is the algorithm used for the approximate projection $\cS_{\bD}$ in the identify and prune steps.

\subsection{Contribution}
As observed in Figure~\ref{sepvclust}, the success of SSCoSaMP depends heavily on both the approximate projection used as well as the structure of the signal; more interesting, none of the methods work well on both signal types.  Unfortunately, this phenomenon is not well understood theoretically and even experimentally.  This motivates the two main points of focus for this paper.  First, we investigate how variants of SSCoSaMP perform on sparse vectors of varying structures.  We rigorously test each individual algorithm on a variety of signal types in order to precisely identify where some algorithms succeed and where others fall short, thus cataloging which methods should be used on which signal classes. 
Secondly, we develop novel CS algorithm variants with the ability of allowing reconstruction over a much wider range of signal structures.  This aims to remove the need to select specific CS methods depending on the structure of the signal of interest.  Our methods improve upon all other existing approaches in this regard.

\section{Empirical Investigation}
The stark contrast in the recovery performance between the SSCoSaMP variants in the clustered and well-separated cases is the main motivation for the rest of this paper. We aim to explore why this sharp rift exists and exactly how far it extends. We reason that the clustered and well-separated cases represent two extremes -- a vast middle ground lay between these two poles that has yet to be explored. To shed some light on this unknown, we create multiple classes of sparse vectors that lay somewhere between these two extremes. In addition to the \textit{clustered} and \textit{spread} signal types of Figure~\ref{sepvclust}, we also define the following classes of sparse signals:

\begin{itemize}
\item \textbf{Hybrid:} a block of $k/2$ nonzeros and an additional $k/2$ nonzeros at least 8 indices away from all other nonzeros
\item \textbf{C Clusters:} $C$ blocks of $k/C$ nonzeros with at least one zero between the blocks
\item \textbf{Alternating:} a block of alternating nonzero and zero entries, for a total length of $2k - 1$ indices
\item \textbf{Pair Spread:} $k/2$ pairs of nonzeros where the nonzeros in each pair are 4 indices apart
\end{itemize}

\textbf{Experimental parameters.} For all of the sparse vectors described above, the actual placement of the nonzero entries is random every trial. We only ensure the the random placement adheres to our general structure. Additionally, in every trial we use a 4$\times$ overcomplete Discrete Fourier Transform (DFT) dictionary with dimensions $256\times 1024$.  The measurement matrix $\bA$ is an $m\times 256$ matrix with standard normal entries.  Unless otherwise stated, the sparsity level is held at $k=8$.  Recovery results for different SSCoSaMP variants can be seen in Table \ref{table:recovery} at the end of Section~\ref{sec:conc}.

\begin{figure}[ht]
\begin{center}
\includegraphics[width=3in]{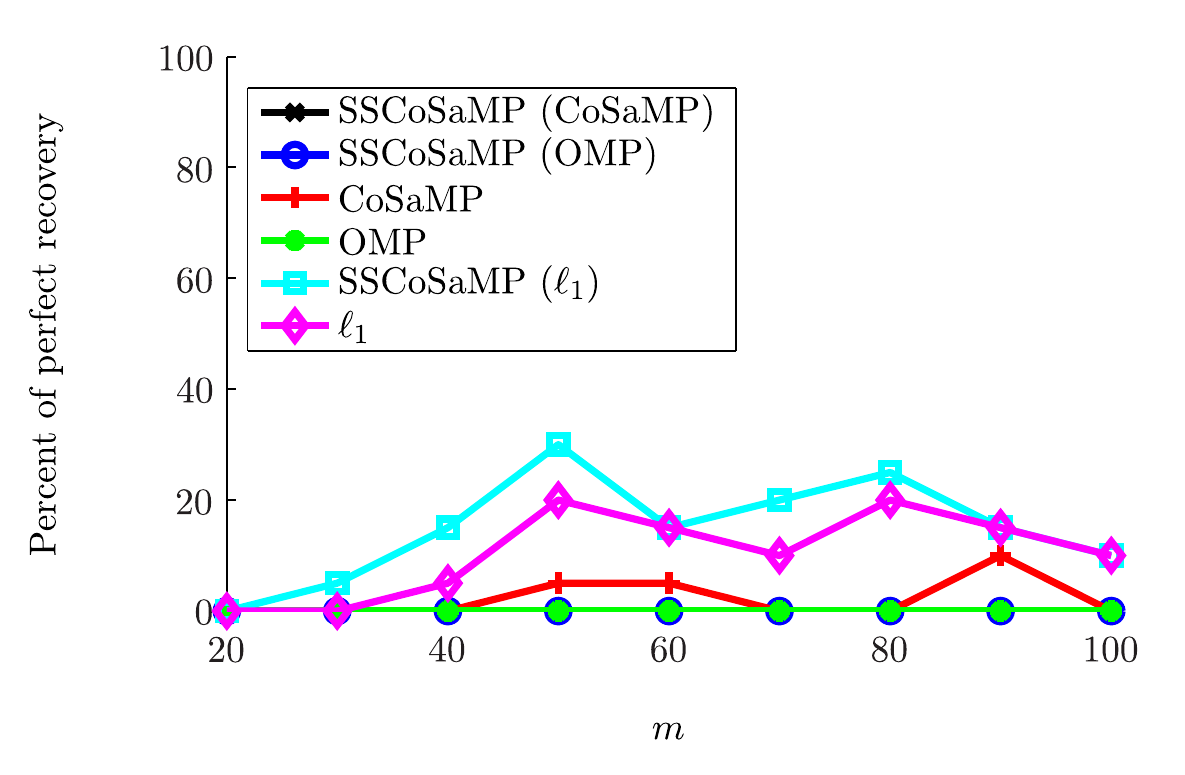}
\caption{Conventional SSCoSaMP and CS algorithms on recovering a sparse vector with a hybrid sparse support: a block of $k/2$ nonzeros with the remaining $k/2$ nonzeros spaced at least 8 slots apart from all other nonzeros.}
\label{fig:hybrid}
\end{center}
\end{figure}

The first, and most telling, sparse vector on which we test the SSCoSaMP algorithms is the hybrid signal, seen in Figure \ref{fig:hybrid}. We see that neither the SSCoSaMP variants nor the classical CS algorithms recover the signal with any reasonable success. While some algorithms do well in a clustered case and others do well with a spread signal, this test shows that none of the algorithms do well when the signal is part clustered and part well-separated. 

\begin{figure}
	\begin{tabular}{cc}
	\includegraphics[width=3in]{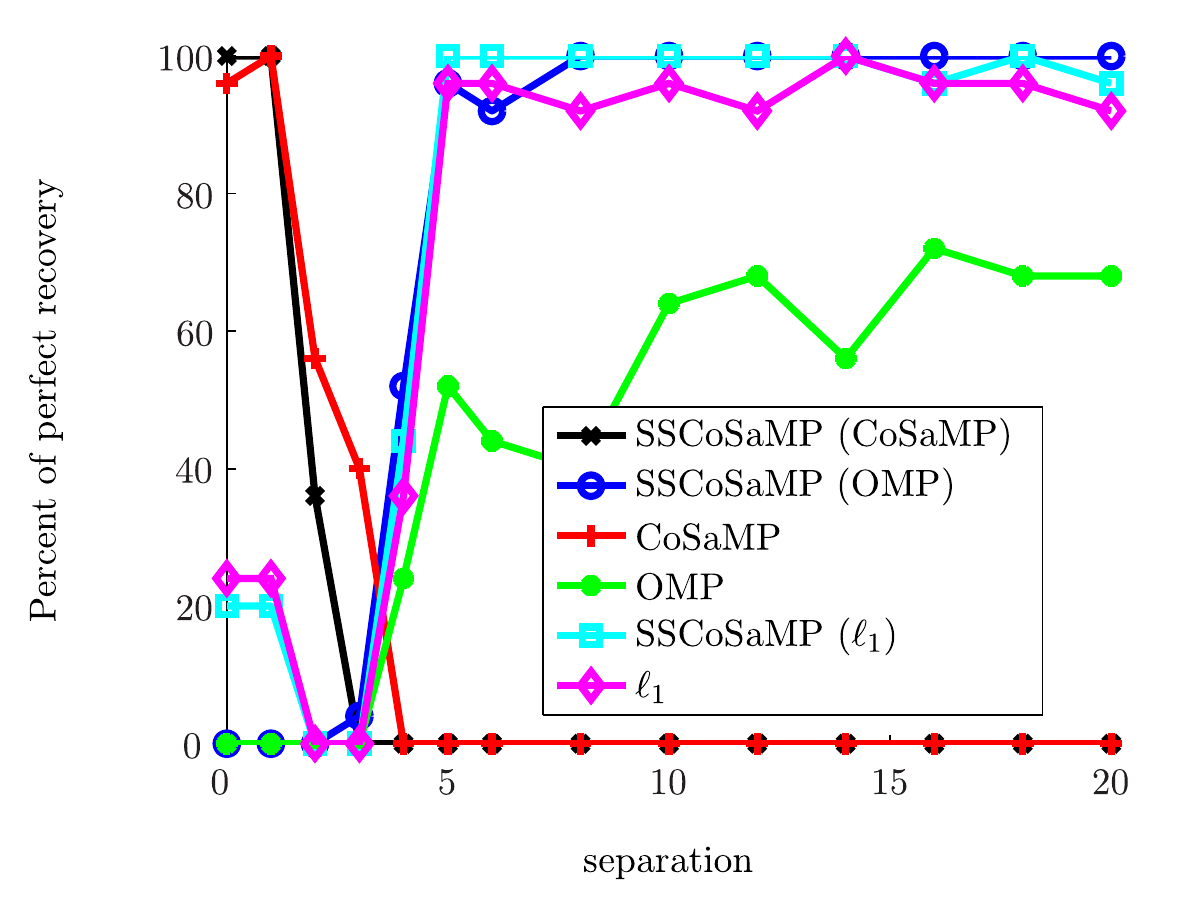} &

	\includegraphics[width=3in]{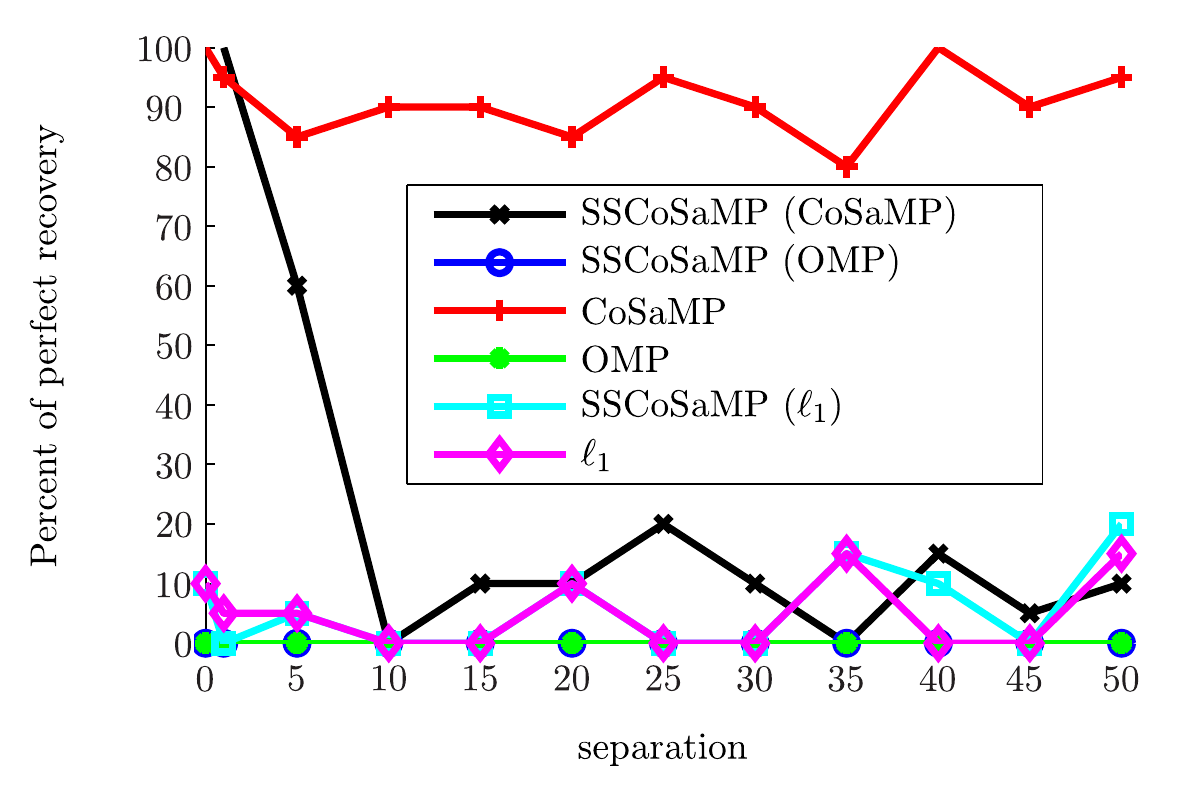}%
\end{tabular}
	\caption {Conventional SSCoSaMP and CS algorithms on recovering a sparse vector. On the left the separations represent the number of zeros between each nonzero entry, and on the right the separations represent the number of zeros between two clusters each of size $k/2$.  Measurements and sparsity are held constant at $m=100$ and $k=8$, respectively. }
	 \label{fig:seps}

\end{figure}

The following test, seen in Figure \ref{fig:seps} (top), summarizes the transition of the performance in the algorithms from a clustered signal to a well-separated one. The $x$-axis now represents the number of zeros between each nonzero entry in a signal. For example, at zero the support of the signal consists of a single cluster of indices, while at a value of 10 separations there are 10 zeros between each nonzero entry. The measurements and sparsity level are held at a constant value, $m=100$ and $k=8$, respectively. 

At a separations value of zero we see SSCoSaMP (CoSaMP) and CoSaMP achieve 100\% perfect recovery, while the other algorithms all perform poorly. As the separations increase to the range of three to four there is a steep drop in the performance of SSCoSaMP (CoSaMP) and CoSaMP. One interesting note is that CoSaMP doesn't drop as fast as SSCoSaMP (CoSaMP) does. 
From five separations onwards SSCoSaMP (CoSaMP) and CoSaMP have no success in recovering the signal. 
When the separations increase to five, SSCoSaMP ($\ell_1$), SSCoSaMP (OMP), and $\ell_1$ start to perform very well. This shows us that in this case the nonzeros need to be at least five spaces apart to be spread enough for these algorithms to recover the signal consistently. We attribute this point to the periodic correlation of columns in the dictionary matrix $\bD$ that is $4\times$ overcomplete. This makes it very difficult for the algorithms to recover the signal at four separations.  Indeed, this required separation is also evident in the super-resolution setting~\cite{Reading2}.  But, as the separations increase past four, the SSCoSaMP ($\ell_1$), SSCoSaMP (OMP), and $\ell_1$ algorithms can recover the signal because they are no longer affected by the periodic correlation of the columns.  

The test seen in Figure \ref{fig:seps} (bottom) is very similar to the previous test. However, the separation value now represents the number of zeros between two clusters of size $k/2$. Again, sparsity and measurements are fixed at $k=8$ and $m=100$, respectively. 
As expected CoSaMP and SSCoSaMP (CoSaMP) perform well for low separation but as the separation between the two clusters increases SSCoSaMP (CoSaMP) starts to perform very poorly. CoSaMP stays around 90\% of perfect recovery for all separation values; it is the only consistent algorithm in this case. SSCoSaMP ($\ell_1$), SSCoSaMP (OMP), OMP, and $\ell_1$ never do well in this test. This result is expected because we have already observed that these algorithms don't perform well when the signal has a clustered type of structure.

\begin{figure}

\centering
\hspace{-.075\linewidth}
\begin{minipage}{.45\linewidth}
\centering
\begin{tabular}{cc}
\raisebox{25mm}{\small \sl {(a)}} &\includegraphics[width=3in]{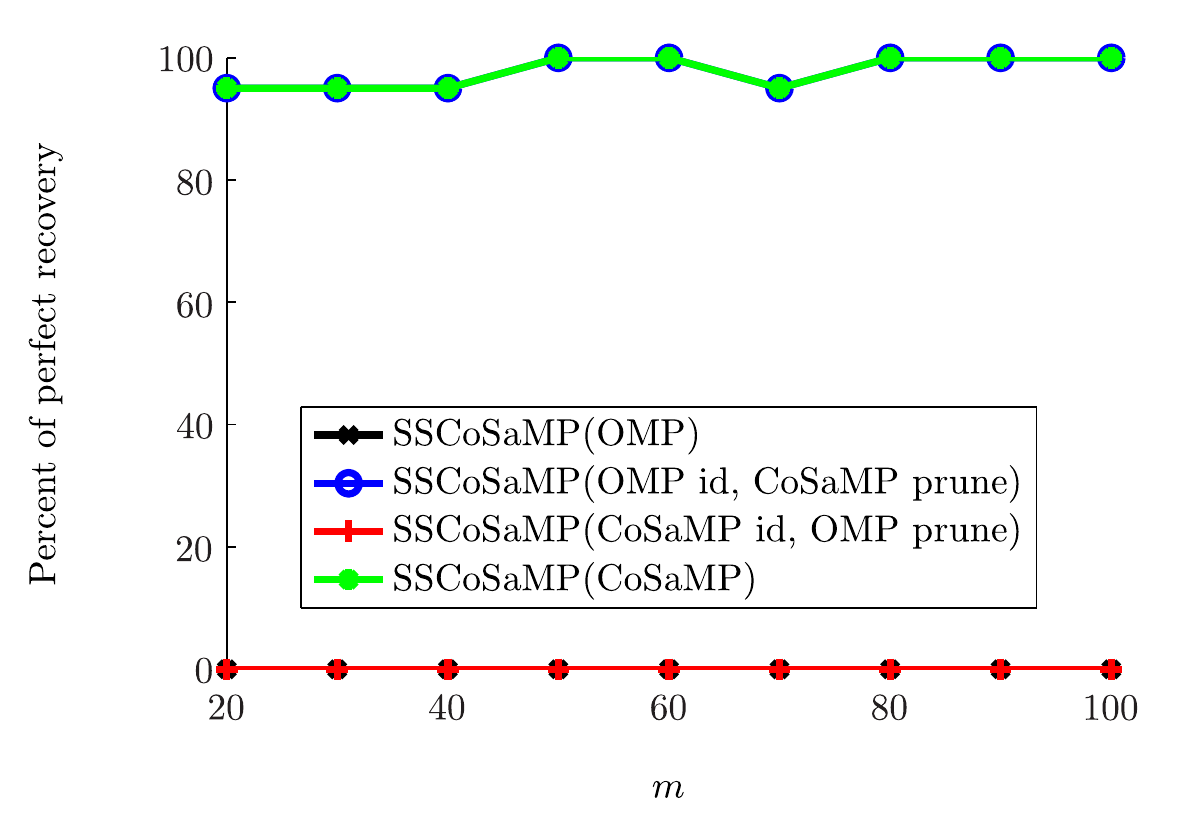}
\end{tabular}
\end{minipage}
\hspace{.05\linewidth}
\begin{minipage}{.45\linewidth}
\centering
\begin{tabular}{cc}
\raisebox{25mm}{\small \sl {(b)}} &\includegraphics[width=3in]{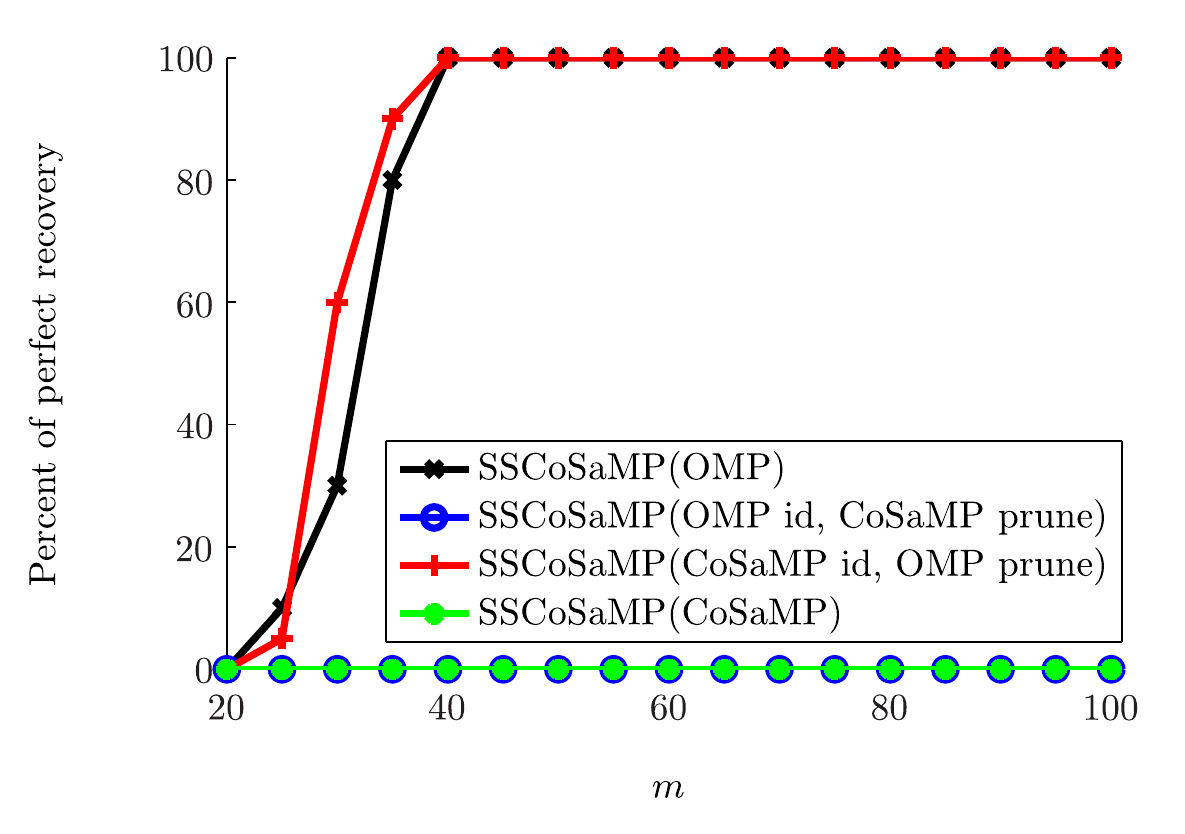}
\end{tabular}
\end{minipage}
   \caption{\small \sl SSCoSaMP variants with differing identification and prune steps, recovery performance on a sparse vector with a single cluster of nonzeros (left) and a sparse vector with well-separated nonzeros (right).
   \label{fig:prunevid}}
\end{figure}

To better understand the limitations of SSCoSaMP, we also examine the inner steps of SSCoSaMP that employ a classical CS method, the prune and identify steps. In Figure \ref{fig:prunevid}, we create two new SSCoSaMP variants. The first uses OMP to identify the near optimal $2k$-sparse support for the proxy vector and CoSaMP to prune the support vector to its best $k$-sparse approximation. The second variant uses CoSaMP to identify and OMP to prune.

We test these two new variants on the two extreme cases of sparse vectors, and the results were definitive. The two new variants behave almost exactly like the pure SSCoSaMP variant whose prune step matches their own. From this, we conclude that the algorithm used to prune determines the overall performance of the algorithm.

\section{Algorithmic Variations}

\subsection{NOMP}
Orthogonal Matching Pursuit (OMP) as described in \cite{OMP}, is a greedy algorithm that uses the signal, $\widetilde{\bv} = \bA^* \bA (\by -\bA\hat{\bx})$ as a proxy for the true desired signal $\bx$.  
OMP constructs the support set by iteratively selecting the coordinate corresponding to the largest coefficient of the proxy $\widetilde{\bv}$. Call this coordinate $\lambda$. OMP iterates $k$ times, where $k$ is the sparsity level. In each iteration, the support and residual are updated. Once the support of the signal $\Gamma$ is found, the signal is estimated using a simple least squares problem as $\hat{\bx} = \bA^\dagger_\Gamma \by$. However, past experiments show that OMP only successfully recovers well-separated signals in the signal space when $\bD$ is an overcomplete dictionary. Here, we present a simple modification in order to recover signals having a variety of structures which are sparse in dictionaries with correlation patterns like the overcomplete DFT. 

The variation of OMP we propose here is deemed Neighborly Orthogonal Matching Pursuit (NOMP). In each iteration of OMP, only one coordinate is added to the support set. In each iteration of NOMP, a $w$-\textit{window} of coordinates is added to the support. That is, a cluster of $w$ adjacent coordinates are added to $\Gamma$. The coordinate corresponding to the largest coefficient, $\lambda$, is at the center of the window.  If $w$ is odd, we select an equal number of adjacent coordinates from either side of the $\lambda$. If $w$ is even, choose the last coordinate by selecting the largest coordinate on the edge of the window. In order to prevent a support set that is larger than the number of measurements $m$, we ensure that $w \leq \frac{m}{k}$. NOMP is described in Algorithm \ref{alg:NOMP}.

NOMP's method of constructing the support set is theoretically similar to the $\epsilon$-OMP method of \cite{epsOMP}. $\epsilon$-OMP assumes the use of a highly coherent dictionary and takes advantage of the correlated columns in $\bD$. $\epsilon$-OMP adds additional coordinates to the support set if these coordinates are correlated with the largest one. That is, $\epsilon$-OMP uses the $\epsilon$-extension defined in \cite{Paper6}, given by:  
\begin{equation} \label{eq:ext}
\mathrm{ext}_{\epsilon, 2} (\Gamma) = \left\{ i : \exists j \in \Gamma, \frac{|\langle d_i, d_j \rangle|^2}{||d_i||^2_2||d_j||^2_2}  \geq 1 - \epsilon^2 \right\}
\end{equation}
This is an alternative way of extending the support set that would have been determined in traditional OMP.  

 When $\bD$ is an overcomplete DFT dictionary, since nearby columns in $\bD$ are correlated, $\epsilon$-OMP essentially adds coordinates neighboring $\lambda$ to the support set. However, it may be difficult to tune the parameter $\epsilon$ in order to select the desired number of nearby neighbors.  For NOMP, one can always reliably set $w = m/k$ if one does not know the structure of the dictionary.

NOMP, however, is designed only for a specific class of dictionaries whose columns have correlated neighbors and is therefore less generalized than $\epsilon$-OMP. NOMP explicitly adds a specified number of additional coordinates adjacent to the largest one. Therefore, we will know the exact size of the resulting support set and can ensure $| \Gamma | \leq m$. If there is knowledge of the approximate size of any clusters in the signal that needs to be recovered, we can change $w$ as needed. 

The main difference between NOMP and $\epsilon$-OMP is in the update step. $\epsilon$-OMP does not use the $\epsilon$-extension when updating the estimated signal. The extension is only used in the least squares step at the end of the algorithm.  NOMP, on the other hand, includes the neighbors in the support in the update step and in the least squares step.  It turns out that this may lead to a difference in performance.

\begin{algorithm}[h]
\caption{Neighborly Orthogonal Matching Pursuit (NOMP)} \label{alg:NOMP}
\begin{algorithmic}
\STATE \textbf{input:} Measurement matrix $\bA$, dictionary $\bD$, measurement vector $\by = \bA \bx$ where $\bx = \bD \balpha$, sparsity level $k$, number of measurements $m$, window size $w\leq m/k$
\STATE \textbf{initialize:} Let support set $\Gamma = \varnothing$, residual $\br = \by$ and counter $\ell = 0$. Let $\bPhi = \bA \bD$. 
\WHILE{$\ell<k$}
\STATE
\begin{tabular}{ll}
\textbf{proxy:} & $\widetilde{\bv}= \bPhi^*\br$\\
\textbf{identify:} & Select largest coordinate $\lambda$ of $\widetilde{\bv}$ and set  \\
			 & $T = \{\lambda\}$. Select the $\frac{w-1}{2}$ coordinates to the   \\
			 & left and right of $\lambda$.  Add these\\
			 &  coordinates to $T$. \\
\textbf{update:} & $\Gamma = \Gamma \cup T$ \\
		& $\widetilde{\balpha} = \argmin_{\bz} ||\by - \bPhi_\Gamma \bz||_2$ \\
		&  $\br = \by - \bPhi \widetilde{\balpha}$   \\
		& $\ell = \ell +1$
\end{tabular}
\ENDWHILE
\STATE 
\begin{tabular}{ll}
\hspace{-3 mm} \textbf{output:} & Recovered coefficent vector, $\hat{\balpha} = \bA^\dagger_\Gamma \by$ \\
				& Recovered signal, $\hat{\bx} = \bD \hat{\balpha}$
\end{tabular}
\end{algorithmic}
\end{algorithm}

\begin{figure}
\begin{tabular}{cc}
    {{\includegraphics[width=3in]{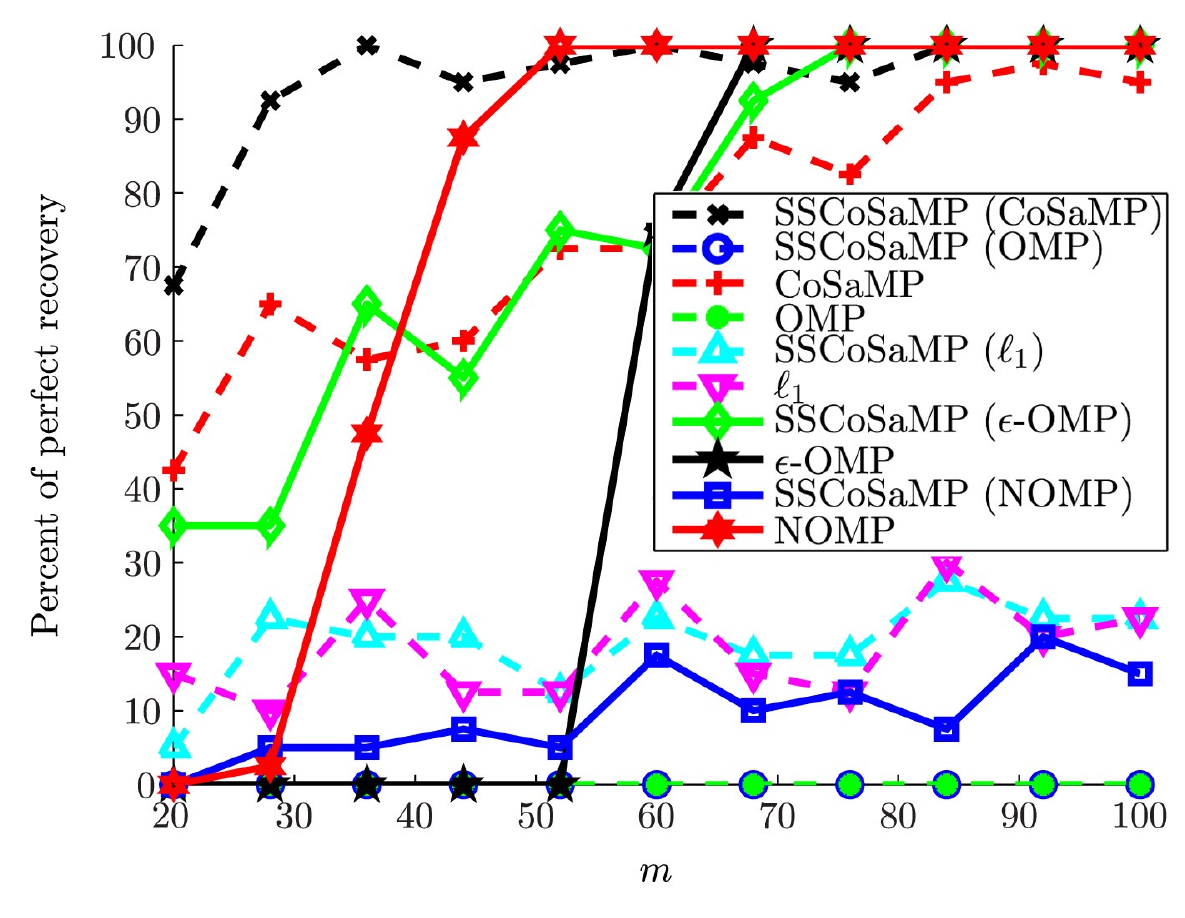} }}&
   {{\includegraphics[width=3in]{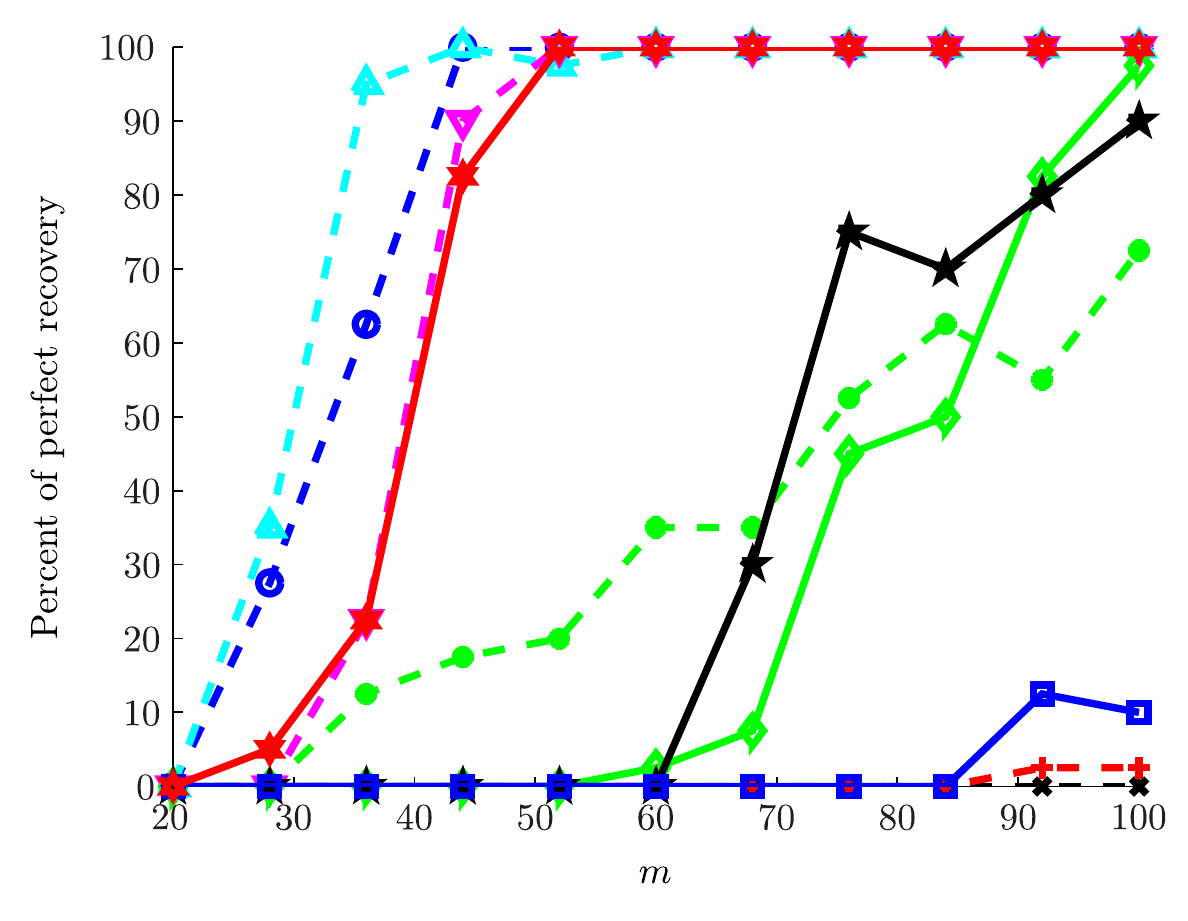} }}\\
    \end{tabular}
    \centering
    \begin{tabular}{c}
   {{\includegraphics[width=3in]{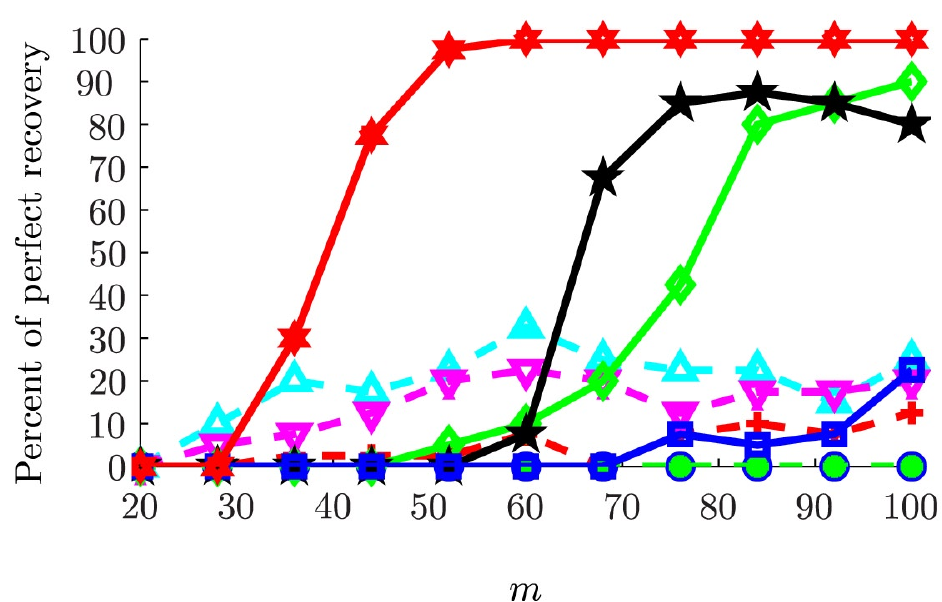} }}\\
     \end{tabular}
	
    \caption{Percent of perfect recovery for clustered signals (left), well-separated signals (right) and hybrid signals (bottom). 
\label{fig:sepclushyb}}
\end{figure}

Our next experiments demonstrate the recovery performance of NOMP benchmarked against $\epsilon$-OMP, traditional compressed sensing algorithms and their SSCoSaMP variants. In all experiments, $k = 8$, $n = 256$, $d = 1024$ and $\bD$ is a $4\times$ overcomplete DFT dictionary. In Figure \ref{fig:sepclushyb}, $m$ is incremented from 20 to 100 by 8 and 40 trials are conducted for each value of $m$. For these tests, the maximum window size is $w = 6$. Figure \ref{fig:sepclushyb} shows experiments on the recovery of clustered signals, well-separated signals and hybrid signals.   NOMP has an excellent performance and is the only algorithm to have $100\%$ perfect recovery for all three signals. Furthermore, it reaches $100\%$ perfect recovery by about $m = 50$ measurements. $\epsilon$-OMP has a decent performance when recovering all three signals. For these tests, $\epsilon = 0.9539$. This value is chosen because when using a $4\times$ overcomplete DFT dictionary, this value corresponds to the selection of six columns. Thus, we can fairly compare NOMP to $\epsilon$-OMP. NOMP's outperformance of $\epsilon$-OMP in these cases highlights the importance of using the neighbors in the update step, and of course $\epsilon$-OMP could also be modified in this way.  We focus on NOMP because of its simplicity in this context.  

\begin{figure}[ht]
\begin{center}
\begin{tabular}{cc}
\includegraphics[width=3in]{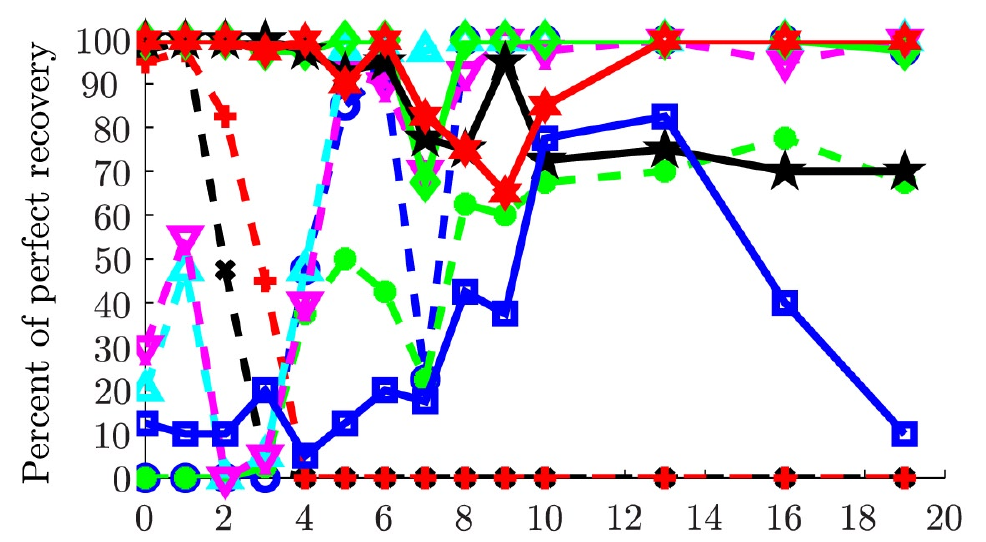} & \includegraphics[width=3in]{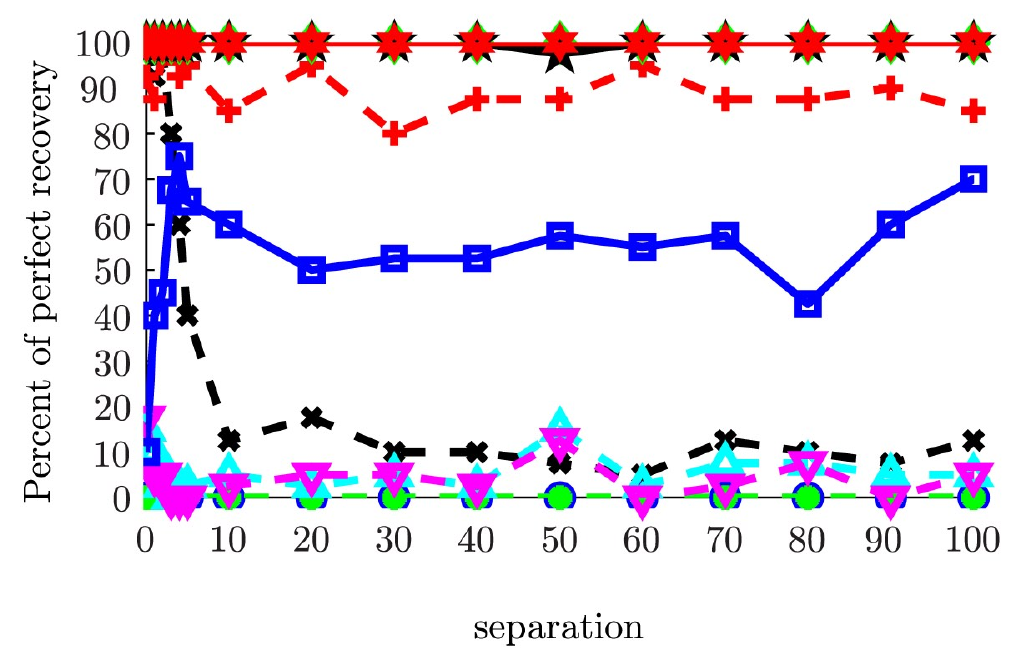}
\end{tabular}
\caption{Left: Percent of perfect recovery of a well-separated signal as a function of the separation between nonzero coefficients. Legend given in Figure~\ref{fig:sepclushyb}. We measure the percent of perfect recovery as the separation between each nonzero coefficient increases.  Right: Percent of perfect recovery of a signal with two clusters as a function of the distance between clusters. Legend given in Figure~\ref{fig:sepclushyb}. We measure the percent of perfect recovery as the separation between the two clusters increases.} 
\label{fig:altsep}
\end{center}
\end{figure}

\begin{figure}[ht]
\begin{center}
\begin{tabular}{cc}
\includegraphics[width=3in]{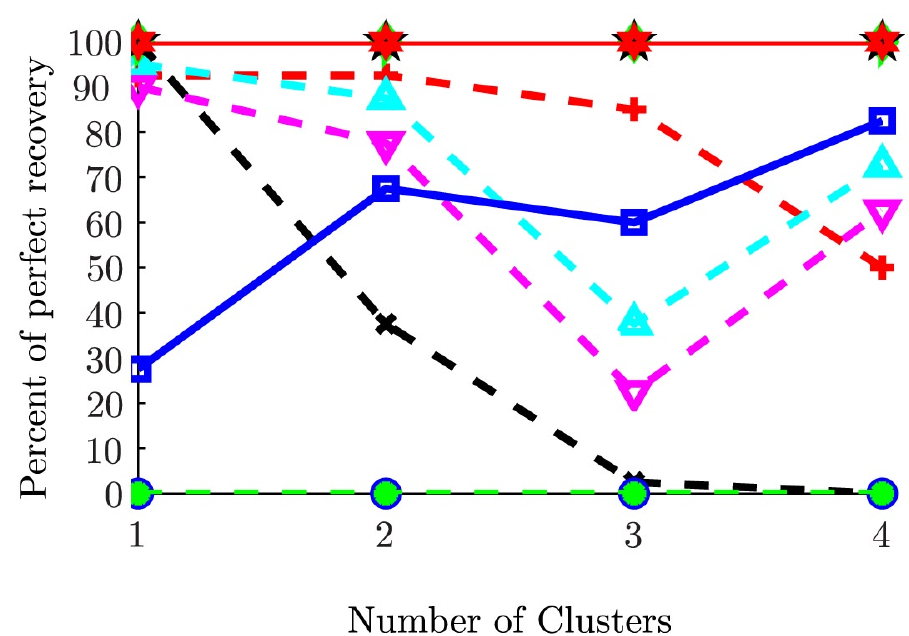} & \includegraphics[width=3.5in]{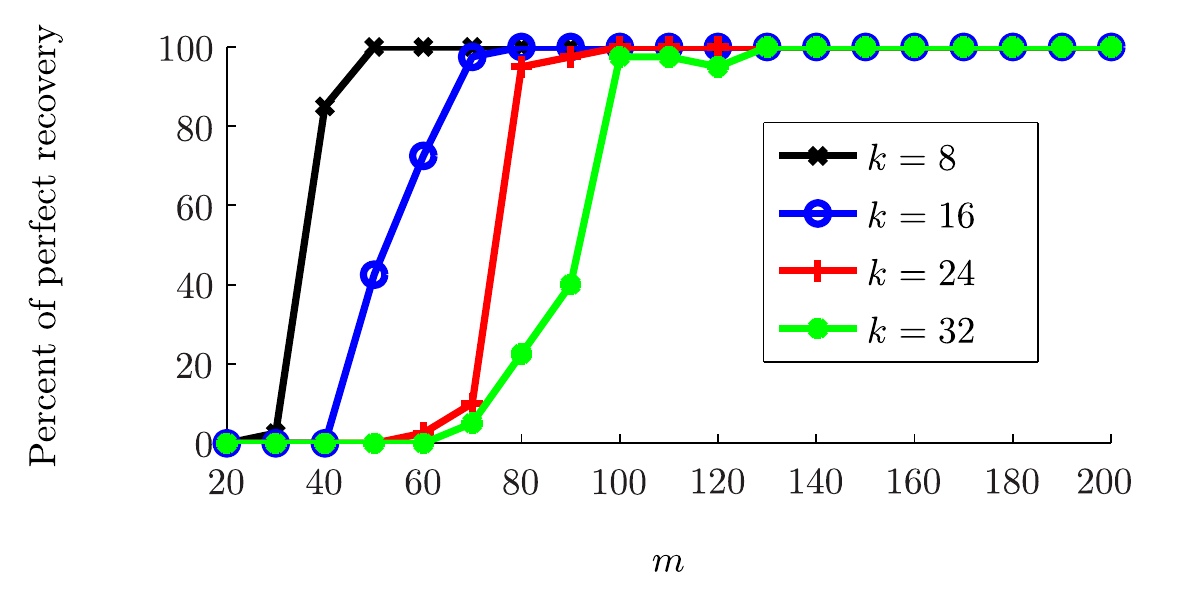}
\end{tabular}
\caption{Left: Percent recovery of a signal with one, two, three and four clusters.  Legend given in Figure~\ref{fig:altsep}.  Right: Percent recovery of NOMP on clustered signals with varying levels of sparsity. }
\label{fig:numclus}
\end{center}
\end{figure}

Figures \ref{fig:altsep} and \ref{fig:numclus} show tests on NOMP and traditional compressed sensing algorithms using a variety of specific signal variants. For all of these tests, the number of measurements is fixed at $m = 100$.  The figure captions provide a description of the type of signal tested. Observe that NOMP and $\epsilon$-OMP have the best performance. Figure \ref{fig:numclus} (right) shows NOMP's success when the sparsity of the signal is increased. 

\begin{figure}
\begin{center}
\includegraphics[width=5in]{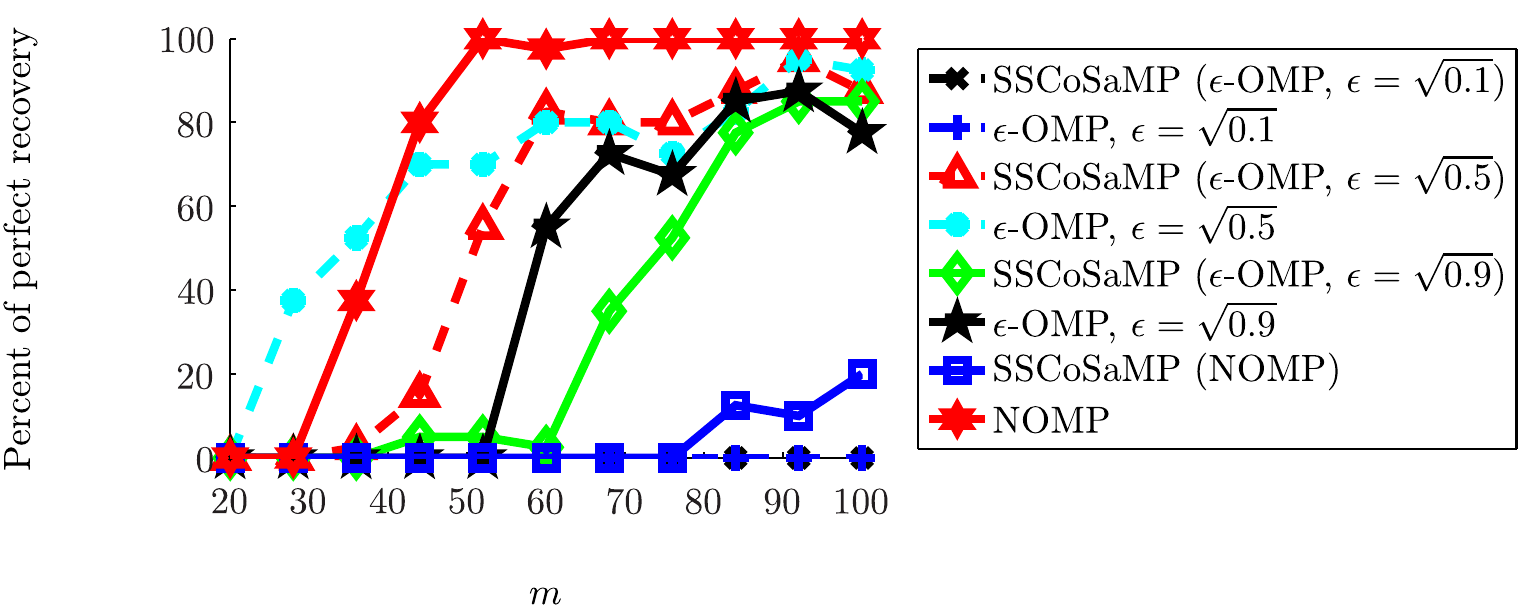}
\caption{Percent of perfect recovery of NOMP and $\epsilon$-OMP for various choices of $\epsilon$ on hybrid signals.  }
\label{fig:nompeps}
\end{center}
\end{figure}

From these experiments, it is clear that NOMP is the algorithm with the best performance when recovering a variety of signal types, which is not surprising since it is designed specifically for this type of dictionary. 
Figure \ref{fig:nompeps} shows the percent of perfect recovery of NOMP and $\epsilon$-OMP with varying values of $\epsilon$. All algorithms are tested on a hybrid signal since Figure  \ref{fig:sepclushyb} shows that NOMP is the only algorithm to successfully recover this type of signal. Figure \ref{fig:nompeps} reveals that it is important to select the correct value of $\epsilon$. However, NOMP has been able to perfectly recover all signals with a window size of $w = 6$. 

Another advantage to NOMP is its fast runtime. Its runtime is comparable to that of OMP because NOMP only alters the way OMP selects the support set. SSCoSaMP is much slower due to its need to run another traditional CS algorithm in its identify and prune steps. 
More importantly, NOMP succeeds on a much wider range of signal structures.  Thus, if the signal structure is unknown or the class of signals of interest includes many types, NOMP offers many practical advantages.

\subsection{USSCoSaMP}

The results of Figure~\ref{sepvclust} suggest that one class of algorithms work well for the clustered signal model and a disjoint set of algorithms work well for the well-separated model. Indeed, when OMP and $\ell_1$-minimization are used to form these projections, SSCoSaMP is able to accurately recover the signal when the nonzero entries in the coefficient vector are well-separated, but unable to recover the signal when the nonzeros are clustered together. Conversely, when CoSaMP is used to form these projections, SSCoSaMP exhibits near-perfect signal recovery when the coefficient vector's nonzeros are clustered, and is unable to recover the signal when the nonzeros are well-separated.

A natural question is whether one can leverage the benefits of each algorithm type simultaneously in order to obtain successful reconstruction of both signal types, as well as those in between.  
For computational efficiency, we choose to focus on the greedy algorithms. After observing the signal recovery success of OMP and CoSaMP in the two ``extreme" cases of coefficient vector structure, we seek to combine their performance to create a method that can succeed at both extremes and everywhere in between.

A natural way to combine the recovery performance of SSCoSaMP (OMP) and SSCoSaMP (CoSaMP) results in the USSCoSaMP algorithm (acrostic: Unionized SSCoSaMP) which is detailed in Algorithm \ref{alg:uss}.

\begin{algorithm}[ht]
\caption{Unionized Signal-Space CoSaMP (USSCoSaMP)}
\begin{algorithmic}

\STATE \textbf{Input:} $\bA$, $\bD$, $\by$, $k$, stopping criterion
\STATE \textbf{Initialize:} $\br = \by$, $\bx_0 = 0$, $\ell = 0$, $\Gamma_{old} = \varnothing$

\WHILE{not converged}
\STATE
\begin{tabular}{ll}

\textbf{Proxy:} & $\widetilde{\bv} = \bA^*\br$ \\
\textbf{Identify:} & $\Omega = \cS_{\bD}(\widetilde{\bv},2k)$\\
\textbf{Merge:} & $T = \Omega \cup \Gamma_{old}$\\
\textbf{Update:} & $\widetilde{\bw} = \argmin_{z}||\by - \bA \bz||_2 \quad   \mathrm{s.t.}  \quad  \bz \in \cR(\bD_T)$\\
\textbf{Prune:} & $\Gamma_{omp} = \cS_{\bD}(\widetilde{\bw}, k)$\\
& $\Gamma_{cosamp} = \cS_{\bD}(\widetilde{\bw}, k)$\\
\textbf{Union:} & $\Gamma = \Gamma_{omp} \cup \Gamma_{cosamp}$\\
\textbf{Update:} & $\widetilde{\bx} = \argmin_{z}||\by - \bA \bz||_2 \quad   \mathrm{s.t.}  \quad  \bz \in \cR(\bD_\Gamma)$\\
 & $\bx_{\ell+1} = \cP_{\Gamma}\widetilde{\bx}$\\
& $\br = \by - \bA\bx_{\ell + 1}$\\
& $\ell = \ell + 1$\\
& $\Gamma_{old} = \Gamma$

\end{tabular}
\ENDWHILE

\STATE \textbf{Output:} $\hat{\bx} = \bx_{\ell}$

\end{algorithmic}
\label{alg:uss}
\end{algorithm}

USSCoSaMP's crucial improvement over SSCoSaMP appears in the prune step. SSCoSaMP makes use of a simple CS algorithm, classically either OMP or CoSaMP, to determine the best $k$-sparse approximation of $\bal$, denoted $\hat{\bal}$, and returns both $\hat{\bal}$ and $\supp(\hat{\bal})$. By contrast, USSCoSaMP uses both OMP and CoSaMP to compute the best $k$-sparse approximation of $\bal$, and only returns the supports of each approximation, denoted in Algorithm~\ref{alg:uss} as $\Gamma_{omp}$ and $\Gamma_{cosamp}$. After taking the union of these supports, $\Gamma$, USSCoSaMP performs the same final step as OMP and CoSaMP by solving a least squares problem restricted to the ``best" columns. Because $k \leq |\Gamma| \leq 2k$, the $\hat{\bal}$ returned by USSCoSaMP is between $k$- and $2k$-sparse. 

Like SSCoSaMP, there are different variants of USSCoSaMP. The key change between these variants lies in the identification step. We focus on two variants that perform the best over the broadest class of sparse coefficient vectors. The first alternates between using OMP and CoSaMP to compute $\Omega$. In this case, $|\Omega| = 2k$ always. The second has OMP and CoSaMP compute $\Omega_{omp}$ and $\Omega_{cosamp}$ respectively, then sets $\Omega = \Omega_{omp} \cup \Omega_{cosamp}$. In this case, $2k \leq |\Omega| \leq 4k$. In the following figures, these variants are denoted USSCoSaMP (alt) and USSCoSaMP (union) respectively.

\begin{figure}[ht]
    \centering
    \begin{tabular}{cc}
    {\includegraphics[width=3in]{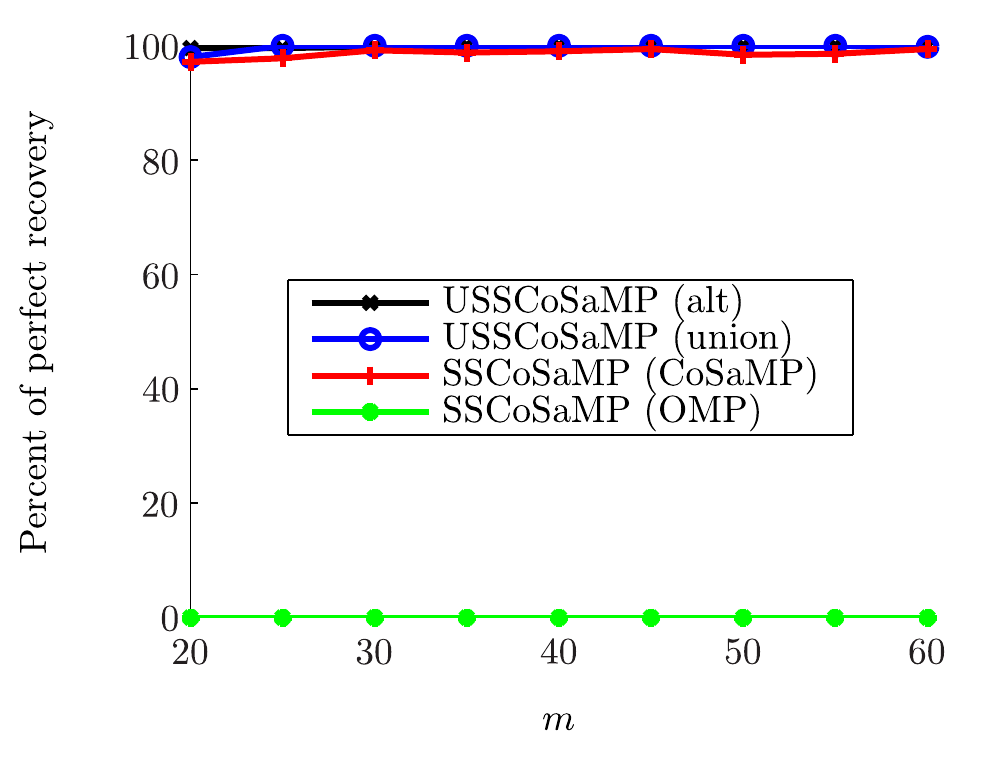} } & 
   {{\includegraphics[width=3in]{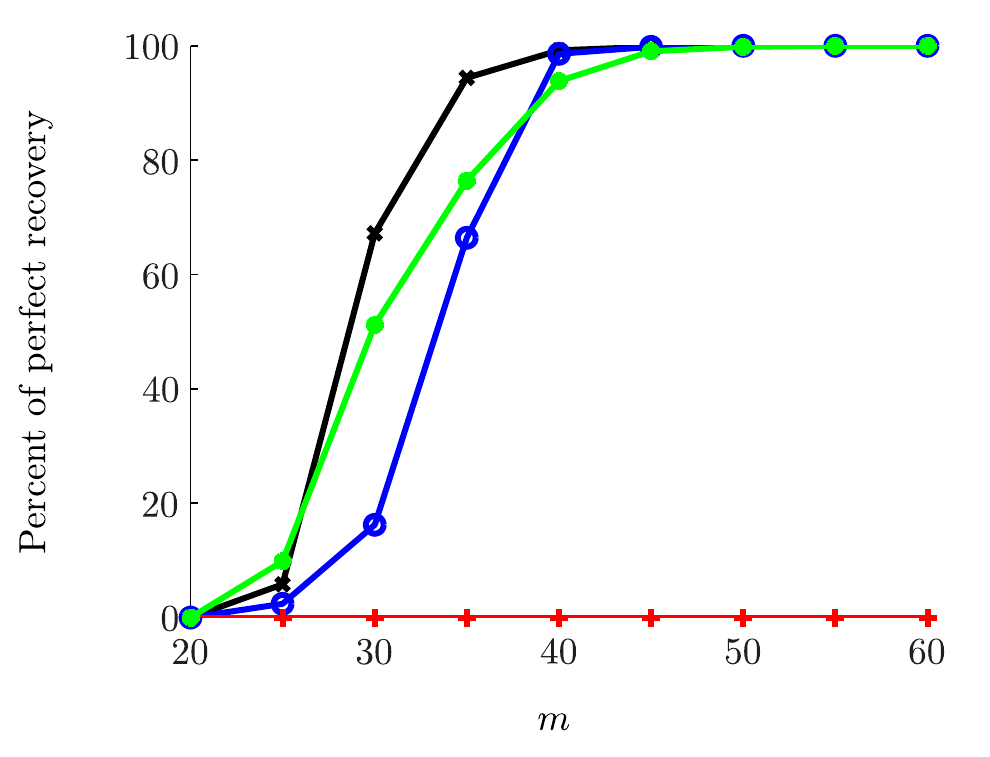} }} \\
   	\end{tabular}
     \centering
    \begin{tabular}{c}
   {{\includegraphics[width=3in]{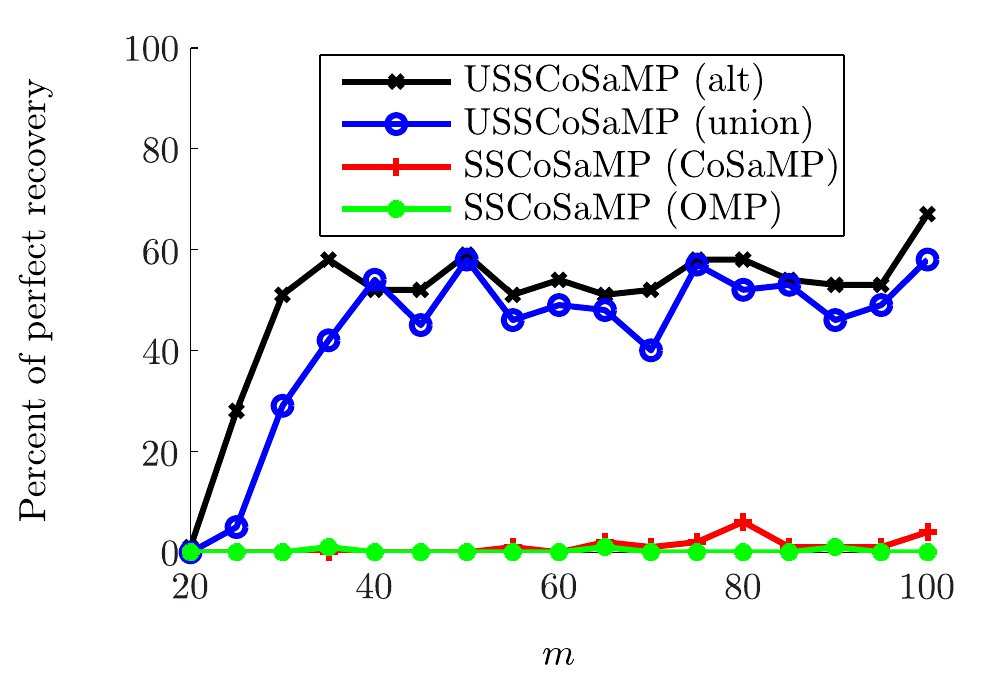} }}
	\end{tabular}
\caption{USSCoSaMP and SSCoSaMP variants recovery performance on sparse clustered signals (left), sparse well-separated signals (right), and hybrid sparse signals (bottom). USSCoSaMP (alt) alternates between using OMP and CoSaMP to compute $\Omega$ whereas USSCoSaMP (union) utilizes both simultaneously and computes $\Omega$ as the union of the two output support sets. 
\label{fig:uss1}}
\end{figure}

Figure~\ref{fig:uss1} shows the recovery performance of USSCoSaMP variants benchmarked against the traditional SSCoSaMP variants for nonzero entries in the coefficient vector that are clustered, well-separated, and hybrid, respectively. For the first two experiments, $k = 8$, $n = 256$, and $d=1024$. The number of measurements $m$ was incremented from 20 to 100 by 5, and 100 trials were conducted for each value of $m$ ($500$ trials for the hybrid signal case). 

Figure~\ref{fig:uss1} (top) shows both USSCoSaMP variants performing as well as the excellent SSCoSaMP (CoSaMP) in the clustered case, boasting near perfect recovery even for low values of $m$. Figure~\ref{fig:uss1} (middle) shows both USSCoSaMP variants performing as well as SSCoSaMP (OMP) for $m \geq 40$ in the well-separated case. In fact, the alternating USSCoSaMP variant even outperforms SSCoSaMP (OMP) in the $m \geq 30$ range, requiring only $m=40$ for $100\%$ recovery versus $m=50$ for SSCoSaMP (OMP).  This latter observation is surprising, since SSCoSaMP (OMP) already performs well on this signal model whereas SSCoSaMP (CoSaMP) performs quite poorly.  Figure~\ref{fig:uss1} (bottom) shows USSCoSaMP sporting a modest recovery performance in the hybrid case, a marked improvement over both SSCoSaMP variants, although still not reaching $100\%$ recovery.  We do note that in the trials where the signal was not accurately recovered, the residual also never reached the lower threshold, so at least one can detect when such a failure has occurred.  This is the only generalized method to our knowledge that performs well on both clustered and separated models simultaneously.  Together, the figures show USSCoSaMP performing better than the sum of its parts.

\section{Conclusion}\label{sec:conc}

In our work, we conduct a rigorous empirical investigation into the recovery performance of $\ell_1$, OMP, CoSaMP and SSCoSaMP on different types of sparse signals. 
We summarize the results of our experiments in Table~\ref{table:recovery}.  The values in the table represent the percent of perfect recovery at $m=100$ measurements. The red values are between 0-39\%, the black represent 40-79\%, and the blue a range of 80-100\%.  We see that the best performers over the broadest class of signals are USSCoSaMP and NOMP. 

\begin{table*}[ht]
  \begin{tabular}{ | l | c | c | c | c | c | c | c | p{2cm} }
    \hline
     
    Method & Clustered & Spread & Hybrid& Two Cluster & Four Cluster & Alternating & Pair Spread\\ \hline \hline
    SSCoSaMP (CoSaMP) &\textcolor{blue}{100}&\textcolor{red}{0} &\textcolor{red}{0}  &\textcolor{red}{20}  &\textcolor{red}{0} &\textcolor{blue}{100} &\textcolor{red} {0}   \\ \hline
    SSCoSaMP ($\ell_1$) &\textcolor{red}{20} &\textcolor{blue}{100}  &\textcolor{red}{30}  &\textcolor{red}{10}  &{60} &\textcolor{red}{25} & \textcolor{red}{35}   \\ \hline
    SSCoSaMP (OMP)& \textcolor{red}{0}&\textcolor{blue}{100}  &\textcolor{red}{0}  &\textcolor{red}{0}  &\textcolor{red}{0} &\textcolor{red}{0} & \textcolor{red}{5} \\ \hline
    CoSaMP&\textcolor{blue}{100} & \textcolor{red}{0} &\textcolor{red}{10}  & \textcolor{blue}{100}  &60 &\textcolor{blue}{100}& \textcolor{red}{0} \\ \hline
    OMP&\textcolor{red}{0} &60  &\textcolor{red}{0}  &\textcolor{red}{0}  &\textcolor{red}{0} &\textcolor{red}{0} & \textcolor{red}{10} \\ \hline
    $\ell_1$&\textcolor{red}{20}  &\textcolor{blue}{100}  &\textcolor{red}{20}  &\textcolor{red}{10}  &50 & \textcolor{red}{20} & \textcolor{red}{25}    \\ \hline
    USSCoSaMP& \textcolor{blue}{100} &\textcolor{blue}{100}  &65 &\textcolor{blue}{80} &\textcolor{red}{30} &\textcolor{blue}{100}&\textcolor{red}{10} \\ \hline
    NOMP&\textcolor{blue}{100}&\textcolor{blue}{100}&\textcolor{blue}{100}&\textcolor{blue}{100} &\textcolor{blue}{100} &\textcolor{blue}{100} &\textcolor{blue}{100}\\ 
    
    \hline
  \end{tabular}
 
  \caption{SSCoSaMP variants and new algorithms' performance on various types of sparse coefficient vectors. All of which are sparse with respect to a $4\times$ overcomplete DFT dictionary. A minimum of 40 trials were performed on each test.}
 \label{table:recovery}  
\end{table*}

The new algorithms that we developed, NOMP and USSCoSaMP, improve upon existing algorithms for recovering signals in a DFT dictionary. Both bridge the gap between clustered and well-seperated signal types and are able to give high recovery on other types of signals as well. Our experimental results provide more insight on the subject and can serve as a guide for practitioners using these methods. NOMP works well for basically all signal types and has a fast runtime, but is specific to certain dictionary structures. USSCoSaMP also provides accurate recovery when the signal is well-seperated or clustered and often on hybrid signals.  We thus provide a catalog of existing and new algorithms along with their performance on different signal structures.  It is important future work to continue the development and analysis of efficient methods that provide accurate recovery of signals in the framework of arbitrary dictionaries. 

\bibliographystyle{plain}

\bibliography{../../../../bib}

\end{document}